\documentclass{amsart}
\RequirePackage{fix-cm}
\usepackage{graphicx}

\usepackage{amsmath}
\usepackage{stmaryrd}
\usepackage{amssymb}

\usepackage[margin=1.2in]{geometry}
\usepackage{multirow}
\usepackage{xcolor}
\newcommand{\nico}[1]{\textcolor{black}{#1}}
\newcommand{\norm}[1]{\left\lVert#1\right\rVert}
\newcommand{\ie}{\textit{i.e.,\ }}
\newcommand{\eg}{\textit{e.g.,\ }}
\newcommand{\dd}{\textrm{d}}

\newcommand{\DOF}{DoF\,}
\newcommand{\LD}{\textrm{LD}}

\usepackage{rotating}

\newtheorem{remark}{Remark}

\usepackage{hyperref}
\begin{document}

\title[Global dynamics visualisation from Lagrangian Descriptors]{
Global dynamics visualisation from Lagrangian Descriptors. 
Applications to discrete and continuous systems.	
}

\author[J.\,Daquin]{
J\'er\^ome Daquin         
}
\address{Department of Mathematics (naXys), $61$ Avenue de Bruxelles, $5000$, Namur, Belgium}
\email{jerome.daquin@unamur.be}

\author[R.\,P\'edenon-Orlanducci]{R\'emi P\'edenon-Orlanducci}
\address{ENSTA Paris, Institut Polytechnique de Paris, $91120$ Palaiseau, France}
\email{remi.pedenon-orlanducci@ensta-paris.fr}

\author[M.\,Agaoglou]{Makrina Agaoglou}
\address{Instituto de Ciencias Matem\'aticas, CSIC, C/Nicol\'as Cabrera $15$, Campus Cantoblanco, $28049$ Madrid, Spain}
\email{makrina.agaoglou@icmat.es}

\author[G.\,Garc\'ia-S\'anchez]{Guillermo Garc\'ia-S\'anchez}
\address{Instituto de Ciencias Matem\'aticas, CSIC, C/Nicol\'as Cabrera $15$, Campus Cantoblanco, $28049$ Madrid, Spain}
\email{guillermo.garcia@icmat.es}

\author[A.M.\,Mancho]{Ana Maria Mancho}
\address{Instituto de Ciencias Matem\'aticas, CSIC, C/Nicol\'as Cabrera $15$, Campus Cantoblanco, $28049$ Madrid, Spain}
\email{a.m.mancho@icmat.es}

\address{}

\date{\today}

\date{\today}

\maketitle

\begin{abstract}
 This paper introduces a new global dynamics and chaos indicator
 based on the method of Lagrangian Descriptor apt for discriminating ordered and deterministic chaotic motions in multidimensional systems. \nico{The selected  implementation of this method} requires only the knowledge of orbits on finite time windows and is free of  
 the computation of the tangent vector dynamics (\ie  variational equations are not needed). 
 To demonstrate its ability in 
 \nico{visualising} different dynamical behaviors, \nico{in particular for highlighting chaotic regions,} several stability maps of classical systems, obtained  with different phase space methods, are reproduced. 
 The benchmark examples are rooted in discrete and continuous nearly-integrable dynamical systems, with prominent features played by resonances. 
 These include the Chirikov standard map, higher dimensional symplectic and volume preserving maps, fundamental models of resonances, and a $3$ degrees-of-freedom nearly-integrable Hamiltonian system with a dense web of resonances. 
 The indicator thus appears to be relevant for understanding phase space transport mediated by resonances in nearly-integrable system, as ubiquitous in celestial mechanics or  astrodynamics.   
\end{abstract}	

\tableofcontents

\section{Introduction}
A modern approach in exploring the phase space structures of a dynamical system is based on the computation of finite time chaos indicator over chosen slices of initial conditions. A central quantity to characterise the
rate at which solutions of initially infinitesimally close initial conditions separate along the dynamics is the largest Lyapunov exponent.  Faster indicators, instead of relying on this time-averaged asymptotic value, exploit the transient period for discriminating the nature of the orbit in a much shorter time \cite{mGu21}. A popular fast variational method is the so-called Fast Lyapunov Indicator (FLI)  \cite{cFr97}. 
The FLI has been extensively used in nearly-integrable settings and orbital dynamics across several astronomical scales, finding applications in the study of the stability of navigation satellites in the near-Earth space environment \cite{jDa22} to the resonant structure of exoplanetary systems \cite{aKy18}, including the demarcation of arches in our Solar System allowing fast transport routes \cite{nT020}. Among the myriad of existing variational indicators, others popular methods include the orthogonal FLI \cite{mFo02,rBa16}, the mean exponential growth of nearby orbits (MEGNO) \cite{pCi16}) or the Generalised Alignment Indices (GALI) \cite{cSk16}.  The interested reader might find additional implementation details and references \nico{with these and other indicators in, \eg \cite{aMo02,cSk10,hSk16,Bib2}}.\\

The Lagrangian Descriptor (LD), rooted in oceanographic studies \cite{jMa09,cMe10}, 
is a relatively recent perspective allowing to reveal phase space structures. 
Alike frequency inspired methods (\eg the frequency analysis method of \cite{jLa93}, or the integrated auto-correlation function \cite{rBa09}), and contrarily to variational methods, the frame does not rely on the dynamics in a vicinity of an initial condition to detect hyperbolicity and chaotic motions. The LD scheme, described with more mathematical rigor in Sect.\,\ref{sec:LD}, builds originally on the so-called $M$-function which evaluates, for a given initial condition, the trajectory length\footnote{Albeit we focus in this contribution on the arc-length LD, let us point out that  several others LDs have been proposed in the literature, \nico{ based on the integral of positive quantities} \nico{along trajectories raised to certain powers \cite{aMa13}, the $p-$norm LD family \cite{cLo15,cLo17},} the action based LD \cite{fMo20}, and more recently the geometrical LD \cite{rPO21} detailled in \ref{app:LDgeo}} computed  over a finite size time-window. 
As a matter of fact, the computation of those arc-lengths, as function of the coordinates, \nico{ is able to locate and reveal the geometrical template organising the phase space. 
	These structures include objects such as separatrices of hyperbolic equilibria, the stable and unstable manifolds of hyperbolic orbits, manifolds of normally hyperbolic manifolds, hyperbolic sets, invariant tori
	or generalisations such as Lagrangian coherent structures \cite{aMa13,cLo15,cLo17,C17,GG18,C19,C19b,mBe20,M21}.}
The heuristic idea driving the ability of the LD to detect hyperbolic objects, as described \nico{in \cite{cMe10,cMe12,aMa13}}, is that trajectories that start and evolve close to each other will have similar arc-lengths, whilst those arc-lengths will change ``abruptly'' when crossing separatrices or other separating objects. \\

Besides its historical roots in Lagrangian transport  in geophysical flows \cite{cMe10}, the LD method 
has been consistently boosted and further developed over the years by the theoretical dynamical chemistry community
\cite{gCr15,gCr16,gCr17,aJu17,mAg19,mAg21}, 
and found further applications in
the study of cardiovascular flows \cite{aDa21} or billiard dynamics \cite{gCa22}
to name but a few.  
Rather surprisingly and as far as we are aware, the method has been left untouched by the celestial mechanics and astrodynamical communities, albeit well versed with resonant, diffusive, manifolds driven or chaotic transport studies \cite{jLa89,aMo95,aMo02,sRo06,cCh18}. 
In this work,  we transfer, deploy and benchmark the LD methodology based on the arc-length metric to nearly-integrable settings where prominent features of the dynamics are shaped by resonances and their possible interactions. 
\nico{\cite{aHa17} have noticed weakness of the dynamical features highlighted by this LD  on geophysical contexts. A diagnosis  that systematically overcomes limitations is presented here.}
The paper is outlined as follows:
\begin{itemize}
	\item In Sect.\,\ref{sec:LD}, we present with mathematical rigor the framework of Lagrangian Descriptors. We discuss aspects related to the regularity of the LD application. 
	We then introduce a second derivative based quantity, denoted $\norm{\Delta \LD}$, proposed as a new global dynamics indicator. 
	\item In Sect.\,\ref{sec:Flow}, we demonstrate the ability of the $\norm{\Delta \LD}$ indicator 
	to reveal the dynamical template of continuous systems prompt to resonate. 
	We consider fundamental models of resonances, such as the integrable pendulum, the $2$ degrees-of-freedom (\DOF) modulated approximation, and a $2$-\DOF \,model eligible to Chirikov's overlap of nearby resonances. Those examples are further completed by demonstrating that the Arnold web of the $3$-\DOF \, Froeschl\'e-Guzzo-Lega model is fully recovered by the $\norm{\Delta \LD}$ analysis. 
	\item In Sect.\,\ref{sec:Map}, we extend our results to discrete nearly-integrable systems. 
	We compute stability maps of the paradigmatic standard map and higher dimensional nearly-integrable symplectic and volume preserving maps. 
	Particular attention is paid to the problem of detecting and portraying the geography and interactions among  resonances. 
	The obtained results are validated through indirect comparisons with the phase space method (\ie the computations of orbits) or to stability maps computed with different variational methods.   
\end{itemize}
We close the paper by summarising our results.



\section{Lagrangian Descriptors and 
	the $\norm{\Delta \LD}$  indicator}\label{sec:LD}

\subsection{Framework of Lagrangian Descriptors}
In the following,  the dynamical systems  considered will  be given by autonomous continuous flows or smooth mappings. 
In the continuous case, our setting will be an $m$-\DOF \, autonomous\footnote{We do not consider explicitly non-autonomous Hamiltonian system, $\dot{x}=J\partial_{x}\mathcal{H}(x,t)$, 
	as for our purpose the time variable $t$ might be treated as an independent variable with a trivial dynamics. 
	In other words, non-autonomous systems are treated as autonomous systems by extending the dimension of the phase space.
} (\ie time independent)
Hamiltonian vector field
reading
\begin{align}\label{eq:EDO}
\dot{x}=J\partial_{x}\mathcal{H}(x), \, x=(p,q) \in D \subset \mathbb{R}^{m} \times \mathbb{R}^{m}, \,
\mathcal{H}: D \to D,
\end{align}
with $\mathcal{H} \in C^{k}, k \ge 1$
and $J$ is the skew symmetric matrix 
\begin{align}
J = 
\begin{pmatrix}
0 & I  \\
-I & 0    
\end{pmatrix},
\end{align}
where $I \in \mathbb{R}^{m,m}$ is the identity matrix. 
In Sect.\,\ref{sec:Flow}, particular emphasis will be given to $m=1$, $m=2$ and $m=3$. Given an initial condition $x_{0} \in D$, the trajectory on the time interval $\mathcal{T}=[t_{0},t]$ is defined as the set  $\{\phi^{\tau}(x_{0})\}_{\tau \in \mathcal{T}}$,  where $\phi^{\tau}$ denotes the flow at time $\tau$ (and supposedly defined on $\mathcal{T}$) associated to the Eq.\,(\ref{eq:EDO}).
The discrete setting deals with mapping taking the form
\begin{align}\label{eq:Mapping}
z_{n+1} = M(z_{n}), \, n \in \mathbb{N},
\end{align}
where $M: S \to S$ is a smooth function, and $S$ is the phase space.
In Sect.\,\ref{sec:Map}, $S$ will have the structure of the cylinder or a product of cylinders. 
Given an initial condition $z_{0}$, the orbit associated to Eq.\,(\ref{eq:Mapping}) is the set of state iterates $\{z_{0},z_{1},z_{2},\cdots\}$. \\

In order to introduce concisely  the notations related to
the LD theory, we focus primarily on the continuous case. 
The concepts and notations are extended to the discrete case in a straightforward way (the reader might find additional details  in \cite{cLo15}). 
For a given $x_{0} \in D$ and a final time $t > 0$, Lagrangian Descriptors take the form
\begin{align}\label{eq:DefLD}
\LD(x_{0},t)=\int_{-t}^{t}\mathcal{G}
\big( 
\dot{x}(\tau)
\big) \, \dd \tau, 
\end{align}
where the choice of the observable $\mathcal{G}$  determines what is averaged along the trajectory speed \cite{aMa13,sBa18}. 
Popular choices in the literature are 
\begin{align}
\mathcal{G}(\dot{x})=\sum_{i=1}^{2m} \vert \dot{x}_{i}\vert^{p}, \, p \in (0,1],
\end{align} 
or
\begin{align}
\mathcal{G}(\dot{x})=
\norm{\dot{x}}_{2}=
\sqrt{\sum_{i=1}^{2m}\dot{x}_{i}^{2}}. 
\end{align}
In the rest of the paper, we adopt the latter choice yielding to    
\begin{align}\label{eq:LD}
\LD(x_{0},t) = \int_{-t}^{t} \norm{\dot{x}(\tau)}_{2}\dd \tau,
\end{align}
which represents the arc-length of the trajectory computed over the time window $[-t,t]$ and passing through $x_{0}$ at time $\tau=0$. \\

When taking a dynamical-systems approach to analyzing the flow, fixed points, periodic orbits, invariant manifolds and their possible stable and unstable manifolds constitute important geometrical backbones of the phase space. 
Let us recall that the stable manifold associated to an hyperbolic equilibria $x_{h}$ of Eq.\,(\ref{eq:EDO}) corresponds to the set
\begin{align}
\mathcal{W}^{s}(x_{h})=
\{
x \in D \, \vert \lim_{\tau \to +\infty}\phi^{\tau}(x)=x_{h}
\},
\end{align}
whilst the unstable manifolds correspond to 
\begin{align}
\mathcal{W}^{u}(x_{h})=
\{
x \in D \, \vert \lim_{\tau \to +\infty}\phi^{-\tau}(x)=x_{h}
\},
\end{align}
with similar definitions for hyperbolic orbits.  
The LDs are able to detect their locations in the phase space and to reconstruct finite pieces of the geometry of those anchors. 
When doing so, it is customary to split Eq.\,(\ref{eq:LD}) into  a
``forward in time'' and 
``backward in time'' contributions as
\begin{align}
\LD(x_{0},t) = \LD^{+}(x_{0},t) + \LD^{-}(x_{0},t),
\end{align}
where
\begin{align}
\LD^{+}(x_{0},t) =  \int_{0}^{t} \norm{\dot{x}(\tau)}_{2}\dd \tau,
\end{align}
and
\begin{align}
\LD^{-}(x_{0},t)=\LD(x_{0},t) - \LD^{+}(x_{0},t).
\end{align}
The computation of $\LD^{+}$ keeps trace of stable manifolds, whilst $\LD^{-}$  highlight unstable manifolds.
When applied to the mapping of Eq.\,(\ref{eq:Mapping}), the discrete analogue of Eq.\,(\ref{eq:LD}) on the time window $\llbracket 0,n \rrbracket$ reads
\begin{align}\label{eq:DefDLD}
\LD(z_{0},n)=
\sum_{j=0}^{n-1}
\sqrt{\sum_{i=1}^{k} (z^{i}_{j+1}-z^{i}_{j})^{2}},
\end{align}
where $z^{i}_{j}$ denotes the $i$-th component of $z=(z^{1},\cdots,z^{k})$ at time $j$.
In the following, we compute the various LDs on time windows $[0,t]$ or $\llbracket 0,n \rrbracket$, for a suitable $t$ and $n$ (confer the discussion on the size of the time window in \ref{app:LDgeo}), simplifying the notation $\LD^{+}$ to simply $\LD$. 

\subsection{Regularity of the LD application}
As clearly stated by \cite{cMe12}, ``the position of the invariant manifolds is not contained on the specific values taken by $M$ but on the positions at which these values change abruptly.'' 
Thus, central to the method of LD is the assessment of its regularity. 
\nico{Besides a few 2D linear or non-linear} models where rigorous results proving 
the loss of regularity on hyperbolic structures 
have been achieved (linear saddle, and rotated version of it, \nico{some non-linear autonomous and nonautonomous settings, or hyperbolic sets  see \cite{aMa13,cLo15,cLo17,GG18})}, and heuristic arguments presented by  \cite{cMe10,aMa13}, no general result regarding the loss of regularity of the LD metric has been established\footnote{
	In more general cases, the claim that the LD metric should be non-differentiable for hyperbolic motions is an ansatz easily observed at the numerical level.}.  
The analytical proofs provided on simple models, whatever the LD formulation used (arc-length, p-norm, or action based), \nico{rely on the linear assumption or, for nonlinear systems, on the existence of appropriate changes of variables under Moser theorem or Hartman-Grobman theorem conditions}. 
In fact, given the explicit knowledge of the flow, the integral (\ref{eq:LD}) can be estimated for large enough time windows and, from this estimation, follows the non-differentiability of the LD function for points belonging to the manifolds \cite{aMa13,cLo17,vGa22,gGa22}.   
The regularity of the LD metric thus entails the possibility to  delineate the hyperbolic structures through the computation of LD fields on chosen slices of initial conditions, and the extraction of the norm of its gradient $\norm{\nabla \LD}$ \nico{\cite{G16,C19,C19b,mBe20}} or higher-order derivatives based diagnostic such as the Sobel or Laplacian filters used in image processing \cite{vGa22,mKa22}. \nico{This paper proposes the  exploitation of this property of   LDs as a simple chaotic indicator. }

\subsection{LDs on integrable Hamiltonian systems and the $\norm{\Delta \LD}$ indicator}
Starting from an integrable Hamiltonian system with $n$-\DOF,  after the introduction of action-angle variables $(I,\phi) \in \mathbb{R}^{n} \times \mathbb{T}^{n}$, the Hamiltonian depends only on the actions and can be generically written as
\begin{align}
\mathcal{H}(I,\phi)=h(I).
\end{align} 
From Hamilton's equations of motion, one derive that the actions are constant, whilst the angles evolve linearly with time at a rate determined by the frequency vector
\begin{align}
\varpi:
\left\{
\begin{aligned}
& \mathbb{R}^{n} \to \mathbb{R}^{n}, \notag \\
& I \mapsto \varpi(I)=\partial_{I}h(I).
\end{aligned}
\right.
\end{align}
Eq.\,(\ref{eq:LD}) on the time window $\mathcal{T}=[0,t]$ becomes 
\begin{align}\label{eq:LDlinearInte}
\LD\big((I_{0},\phi_{0}),t\big) = \sqrt{\varpi_{1}(I_{0})^{2}
	+ \cdots +
	\varpi_{n}(I_{0})^{2}}\, \,t.
\end{align}	
For the particular case $n=1$, the LD grows linearly with time at a rate depending on the frequency vector as
\begin{align}
\LD\big((I_{0},\phi_{0}),t\big) = \varpi_{1}(I_{0})
\,t.
\end{align}	  
Leveraging further on this estimation, in order to quantify the regularity of the LD metric
to assess the chaoticity of the orbits, we are led to
estimate numerically both the existence and the magnitude of the
first few derivatives. As it will be clear in the subsequent (confer
remark \ref{rmk:O1} below), the first derivatives of the LD metric might sill
be $\mathcal{O}(1)$,
 and in order to balance the linear trend of its growth, we find convenient to introduce a scalar diagnostic based on the second derivatives of the LD
Let us denote by $x=(x_{1},\dots,x_{n})$ the initial condition, and let $t > 0$ be the final time. One introduces the scalar  
\begin{align}\label{eq:NormLap}
\norm{\Delta \LD(x,t)} = 
\sum_{i=1}^{n} 
\left|
\frac{\partial^{2} \LD(x,t)}{\partial x_{i}^{2}}
\right|.
\end{align}
This scalar is useful in quantifying the regularity of the LD metric, and defines a new global chaos indicator as it will be demonstrated in the rest of the paper on a series of dynamical models.  

\begin{remark}
	Assume the LDs have been computed on a regular discretised 1-dimensional section $\Sigma=[a,b] \subset \mathbb{R}$. 
	The discretised  points $\{\sigma_{j}\}_{j=0}^{N-1}$ of the mesh 
	are given by $\sigma_{j+1}=\sigma_{j}+ h$, $j=0,\dots,N-1$, $h=(b-a)/N$, $\sigma_{0}=a$, $\sigma_{N}=b$. 
	From the set of points $\{\LD(\sigma_{j})\}_{j=1}^{N-1}$, Eq.\,(\ref{eq:NormLap}) is estimated 
	using the second symmetric derivative formula, reading
	\begin{align}
	h^{2}\LD\,''(\sigma_{j})
	\simeq 
	\LD(\sigma_{j+1})+\LD(\sigma_{j-1})-2LD(\sigma_{j}).
	\end{align}
	For the boundary points $a=\sigma_{0}$ and $b=\sigma_{N}$, one uses respectively the formula 
	\begin{align}
	h^{2}\LD\,''(a)
	\simeq 
	\LD(a)-2\LD(\sigma_{1})+\LD(\sigma_{2}),
	\end{align}
	or 
	\begin{align}
	h^{2}\LD\,''(b)
	\simeq 
	\LD(b)-2\LD(\sigma_{N-1})+\LD(\sigma_{N-2}).
	\end{align}
	This approach can be extended to higher dimensional LD fields. Note that, in order to compute $\norm{\Delta \LD}$, we do not resample the initial mesh of initial conditions.
	We underline that the measure of regularity based on the second derivatives is also central to the frequency analysis method  \cite{jLa93}.  In the following, when dealing with 2D sections, we assume the
	resolution $h_{x}$ and $h_{y}$ in the $x-y$ direction respectively to be identical, $h_{x}=h_{y}=h$. 
	Up to an offset in the final value of the index
	$\norm{\Delta \LD}$, we consider generically $h=1$.  
\end{remark}

\begin{remark}\label{rmk:O1}
	As it will be further exemplified, \eg in \nico{Figs.\,\ref{fig:fig1}, \ref{fig:fig7} } and Sect.\,\ref{sub:GFM}, 
	the first derivatives of the LDs on the probed models, on regular domains, might still be $\mathcal{O}(1)$ and thus inappropriate to reveal the global dynamics of the problem through a heatmap. In this respect, Eq.\,(\ref{eq:DefDLD}) is better suited to restore sharply the  geometrical template organising the dynamics, even though the truncation error between the two-point central difference and the second symmetric derivatives formula are both $\mathcal{O}(h^{2})$.   
\end{remark}

\section{Applications to flows}\label{sec:Flow}
This section demonstrates the ability of the $\norm{\Delta \LD}$ indicator to reveal the phase space structures of continuous models. 
We illustrate this on archetypal resonant problems, ranging from the integrable pendulum problem to a model eligible to Chirikov's overlap, including the modulated pendulum approximation. We then demonstrate the ability of the indicator to reveal the geography of the resonance on   
on the $3$-\DOF \, Froeschl\'e-Guzzo-Lega Hamiltonian. 
Applications to the well-studied  $2$-\DOF \,H\'enon-Heiles system are presented in \ref{app:HH}

\subsection{Fundamental models of resonances}\label{subsec:pend}
We first illustrate the driving principles of the LDs using fundamental models of resonances, such as the integrable pendulum model (the first fundamental model of resonance \cite{sBr03}) and higher dimensional complication of it supporting chaotic motions. 
Those models, albeit being relatively ``simple,'' contain the fundamental germs driving the LD metric and illustrating the needs for the second-derivatives based indicator that we introduced in Eq.\,(\ref{eq:DefDLD}).  
We consider next the following three Hamiltonians 
\begin{align}
\left\{
\begin{aligned}
&\mathcal{H}(I,\phi)=\frac{I^{2}}{2}-\cos \phi, \, (I,\phi) \in \mathcal{C}, \, \mathcal{C}=\mathbb{R} \times [0,2\pi], \notag \\ 
&\mathcal{H}_{\mu}(I,\phi,t)=\frac{I^{2}}{2}- \big(1 + \mu \sin t\big)\cos \phi, \, (I,\phi,t) \in \mathcal{C} \times \mathbb{R}, \notag \\
&\mathcal{K}_{\epsilon,\mu}(I,\phi,t)
=
\frac{I^{2}}{2}-\frac{I^{3}}{3}-\frac{\epsilon}{12}
\cos \phi + \mu \cos(2I+\phi+t), (I,\phi,t) \in \mathcal{C} \times \mathbb{R},
\end{aligned}
\right.
\end{align}
where $\mu$ and $\epsilon$  are real parameters.
The Hamiltonian $\mathcal{H}$ is the $1$-\DOF \, integrable pendulum model. 
The Hamiltonian $\mathcal{H}_{\mu}$ is the modulated pendulum,  
corresponding to a pendulum with a frequency varying periodically. This model is paradigmatic for resonances having overlapped completely \cite{aMo02}. 
When $\mu=0$, $\mathcal{H}_{\mu}$ reduces to $\mathcal{H}$. 
The Hamiltonian $\mathcal{K}_{\epsilon,\mu}$ is taken from \cite{jFe18} (unpublished work) and represents a perturbation of the integrable Hamiltonian $\mathcal{K}_{\epsilon,0}$, 
eligible to Chirikov's overlap criterion \cite{bCh79}.
For each model, we compute and discuss the properties of the LD metric  over slices of initial conditions. In the following, the time window is $\mathcal{T}=[0,100]$.  

\subsubsection{LDs on the integrable pendulum $\mathcal{H}$}\label{subsec:LDpend} 
The phase space of the pendulum   using the level set method is shown in the top left panel of Fig.\,\ref{fig:fig1}. The phase space  contains an elliptic fixed point at $(0,0)$ and one hyperbolic saddle at $(\pi,0)=(-\pi,0)$. The separatrix, \ie the energy curve associated  to the hyperbolic equilibria, separates the phase space in motions with distinct qualitative features.  
The cat-eye is filled with librational curves, whilst  outside of it the phase space is foliated  by circulational curves, enclosing the cylinder. 
The resonance aperture, \ie the distance between $I=0$ and the apex of the separatrix, has a width $\delta I$ satisfying 
\begin{align}
\mathcal{H}(\delta I,0)=-1.
\end{align}
Solving the last equation for $\delta I$, one finds $\delta I=2$, \ie the full width of the cat-eye is 
$\Delta I = 2 \delta I =4$. 
The computation of the LDs along a resolved line of initial conditions given by
$\phi=0$ and $I \in [-2.5,2.5]$ (blue dashed line in the phase space of the pendulum)   is shown in the middle row of the top line of Fig.\,\ref{fig:fig1}. 
The landscape contains the fingerprint of the symmetry of the Hamiltonian, $\mathcal{H}(I,\phi)=\mathcal{H}(-I,\phi)$, $I \ge 0$. 
The graph of the LD is regular with respect to $I$, except  at $I=0$ and $I=\pm\delta I =\pm 2$ which correspond respectively to the location of the stable equilibrium and the separatrix crossing. 
The graph of the LD metric allows a precise numerical estimate of the resonance width. The last panel of Fig.\,\ref{fig:fig1} shows the heatmap of the LD field computed over a regular $500\times500$ cartesian mesh of initial conditions. The structures of the phase space are recognised by the LD field. 

\subsubsection{LDs on the modulated pendulum $\mathcal{H}_{\mu}$} 
The phase space obtained through iterations of the period-map (snapshots of the flow  at every multiple period of time $T=2\pi$), the LD metric and the heatmap of the LD field are shown in the second row of Fig.\,\ref{fig:fig1}.  
The computations refer to  $\mathcal{H}_{\mu}$ with $\mu=0.1$ in the extended phase space, \ie to   
\begin{align}\label{eq:WeakPend2DOF}
\mathcal{H}_{\mu}(I,J,\phi,\tau)=
\frac{I^{2}}{2} + J - \big(1 + \mu \sin \tau\big)\cos \phi.
\end{align}
The $\phi-I$ plane (with $J=\tau=0$) contains elliptic and hyperbolic periodic orbits. The unstable periodic orbit  generates the chaotic layer distributed around the unperturbed separatrix of $\mathcal{H}_{0}$.  
The LD metric computed over the line $\phi=0$ is smooth, and becomes irregular when crossing the chaotic layer. 
A notable difference with the integrable case is that the periodic orbit at the origin (the elliptic equilibrium for $\mathcal{H}_{0}$) is no longer a cusp point of the LD map.  The LD field identifies the elliptic region, but fails to reveal sharply the chaotic layer. 

\subsubsection{LDs on the models $\mathcal{K}_{\epsilon,\mu}$ eligible to Chirikov's overlap} 
\begin{figure}
	\centering
	\includegraphics[width=1\linewidth]{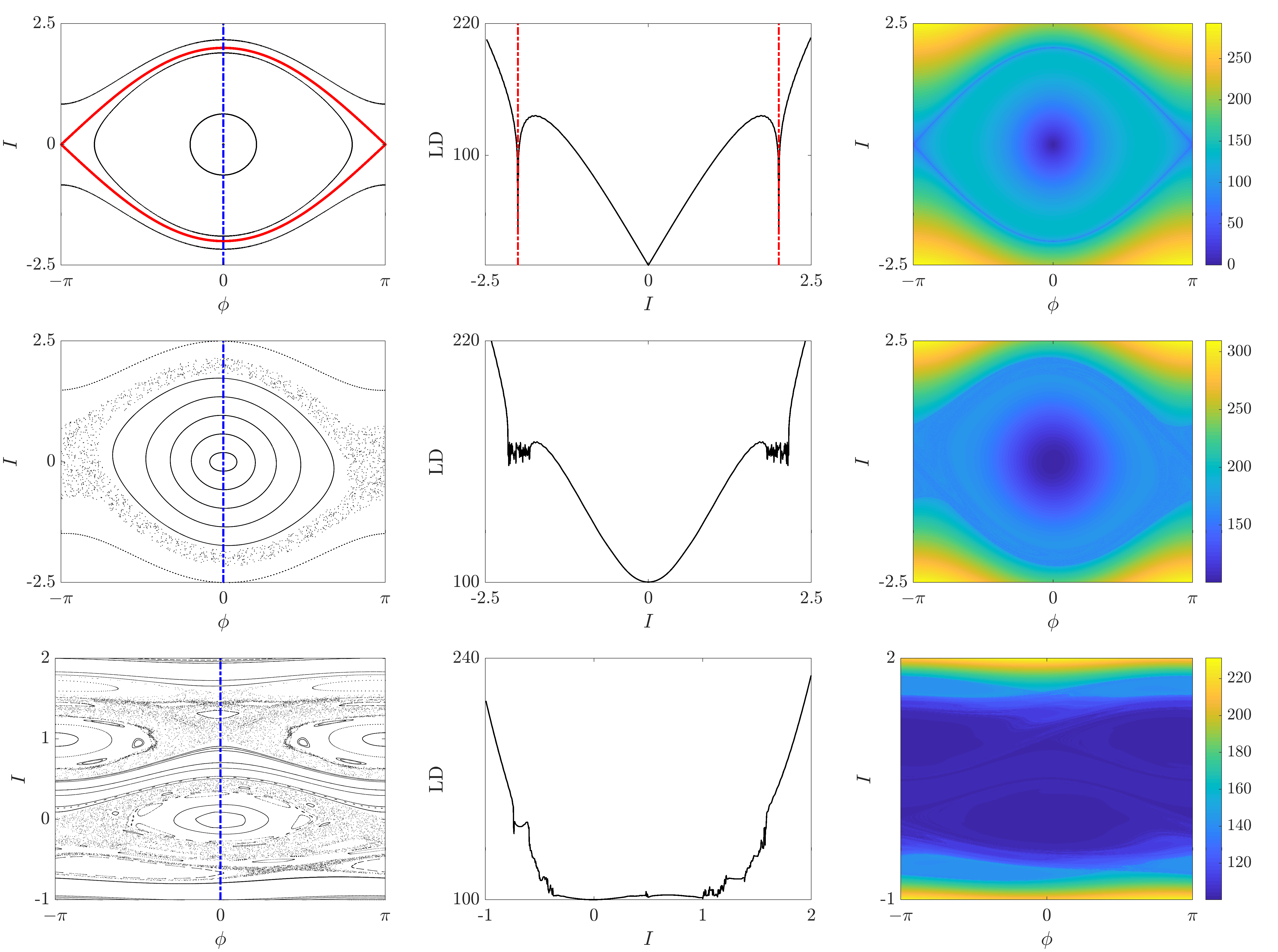}
	\caption{
		The panel shows the phase space (obtained either through the level-set method or iterations of the period-map), the LD landscape computed over the blue dashed line of initial conditions and the LD fields for respectively, from top to bottom, the pendulum model $\mathcal{H}$, the pendulum with varying length $\mathcal{H}_{\mu}$, $\mu=0.1$, and the Hamiltonian $\mathcal{K}_{\epsilon,\mu}$, $\epsilon=0.5$, $\mu=0.01$ eligible to Chirikov's overlap. Due to the overall linear trend of the LD metric with respect to the initial action, the local irregularities of the LD landscape when crossing hyperbolic domains  are hard to detect. As a result, in case of a rich dynamical template, the heatmap of the LD field is ``flat'' and does not reveal sharply the dynamical structures.
	}\label{fig:fig1}
\end{figure}

We now proceed in a similar way to the analysis of  $\mathcal{K}_{\epsilon,\mu}$, $\epsilon=0.5$, $\mu=0.01$. The information are gathered  in the third row of Fig.\,\ref{fig:fig1}.  The integrable dynamics $\mathcal{K}_{\epsilon,0}$ contains resonant eyes filled with librational curves centered around the resonant actions $I=0$ and $I=1$. The widths of the resonant islands are comparable to their mutual distance. For $\mu \neq 0$, resonances interact leading to the apparition of chaotic motions and secondary resonances. 
From the LD landscape regularity, one still guess the location of the hyperbolic structures. Nevertheless,  
this landscape contains two different scales making the appreciation hard: on one side, there is the local and confined loss of regularity and, on the other side, the LDs values spreading a large domain.   
As a result, the heatmap of the LD field is unstructured, and  fail in restoring  the dynamical template offered by the period-map.  

\begin{figure}
	\centering
	\includegraphics[width=0.85\linewidth]{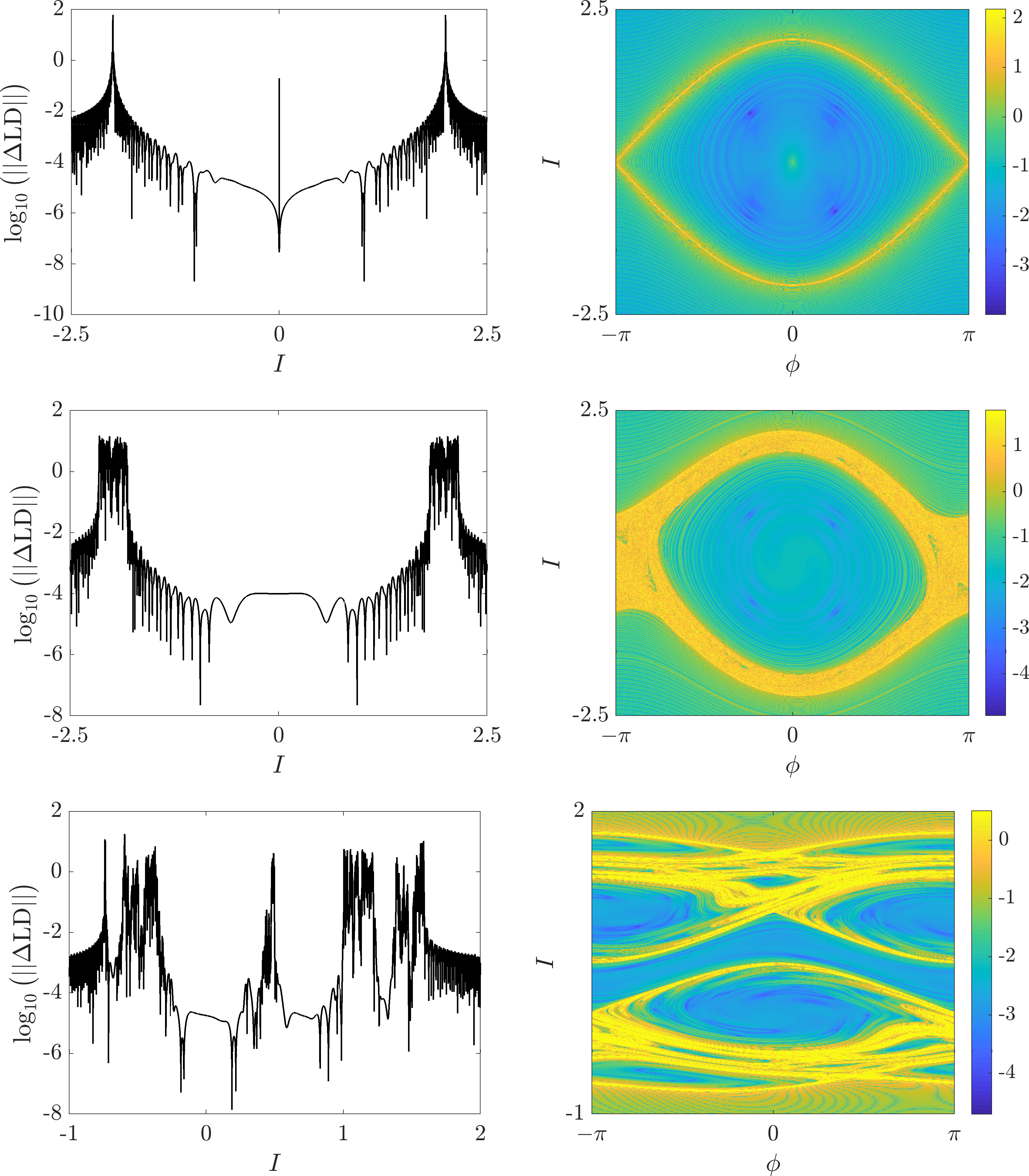}
	\caption{
		Same as Fig.\,\ref{fig:fig1} using 	 $\log_{10}\big(\norm{\Delta \LD}\big)$ instead of LD. Contrarily to the LD metric, the $\norm{\Delta \LD}$ indicator delineates sharply the dynamical structures of the phase space. 
	}\label{fig:fig2}
\end{figure}

Fig.\,\ref{fig:fig2} is the analogue of Fig.\,\ref{fig:fig1}, keeping the numerical settings unchanged, but using $\log_{10}(\norm{\Delta \LD})$ instead of LD. 
The passage from LD to $\norm{\Delta \LD}$ landscapes annihilates the linear trend in the LD metric and clearly emphasises the location of the hyperbolic structures  or chaotic domains, by taking values different by several orders of magnitude compared to the values taken on regular motions (except for elliptic equilibria, where the LD metric is also irregular). As a result, the $\norm{\Delta \LD}$ indicator reinflates the former LD fields, and the formerly missing dynamical structures are now clearly distinguishable. 

\begin{remark}[Oscillations induced by the derivatives]
	As one observes in Fig.\,\ref{fig:fig2}, the computation of the second derivatives introduces oscillations in the $\norm{\Delta \LD}$ landscape (albeit absent in the LD landscape) and, further, Moir\'e-like patterns in the $\norm{\Delta \LD}$ fields. 
	Those oscillations appear primarily within the regular domains, where the derivatives oscillate 
	by about two orders of magnitude. 
	On the other hand, we have computed the same landscapes on  much more resolved grids of initial conditions and we have observed that the values taken by $\norm{\Delta \LD}$ decrease significantly within the regular domains. Altogether, the structures have no dynamical significance  and are of least importance for our goal of highlighting hyperbolic and chaotic domains.
\end{remark}

In the remaining sections, we provide further evidences that the $\norm{\Delta \LD}$ indicator succeeds in portraying the fine distribution of ordered and chaotic motions on other well-studied examples of the literature.

\subsection{The Froeschl\'e-Guzzo-Lega Hamiltonian}\label{subsec:FGL}
\begin{figure}
	\centering
	\includegraphics[width=1\linewidth]{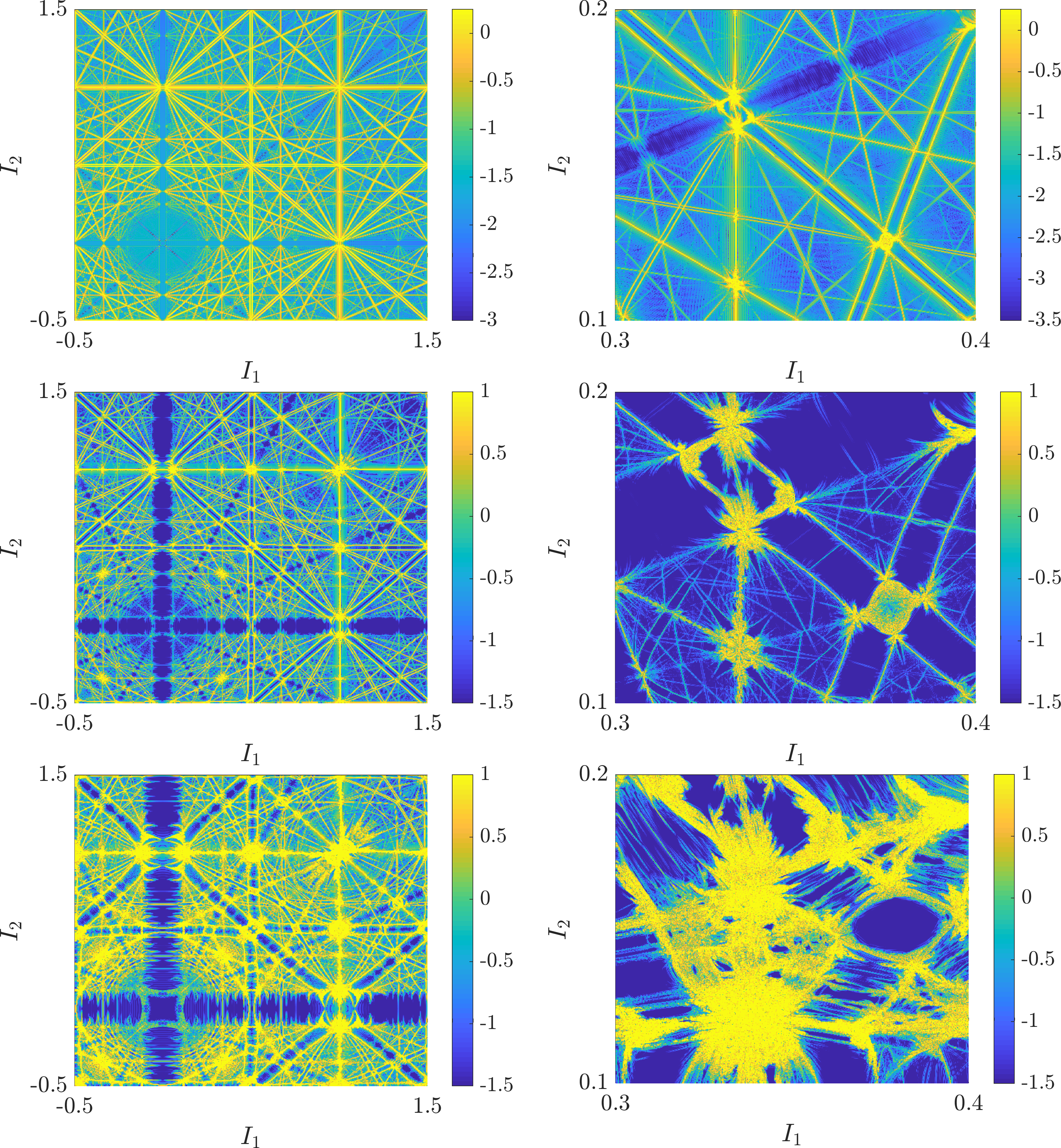}
	\caption{
		$\norm{\Delta \LD}$ stability maps associated to the 
		Froeschl\'e-Guzzo-Lega Hamiltonian of 
		Eq.\,(\ref{eq:FGL}) for 
		$\epsilon=0.001$ (top line), 
		$\epsilon=0.01$ (middle line) and
		$\epsilon=0.04$ (bottom line). Each plot in the right-hand side is a zoomed-in portion of the phase space explored on the left-hand side. 
	}\label{fig:fig3}
\end{figure}

The Froeschl\'e-Guzzo-Lega Hamiltonian corresponds to the $3$-\DOF \, Hamiltonian function
\begin{align}\label{eq:FGL}
\mathcal{H}_{\epsilon}(I_{1},I_{2},I_{3},\phi_{1},\phi_{2},\phi_{3})=
\frac{I_{1}^{2}}{2}+\frac{I_{2}^{2}}{2}+I_{3}+
\frac{\epsilon}{(\cos \phi_{1} + \cos \phi_{2} + \cos \phi_{3} +4)}, 
\end{align}
where $(I,\phi) \in \mathbb{R}^{3} \times \mathbb{T}^{3}$ and $\epsilon \in \mathbb{R}$ is a parameter.  
This system has been studied in a number of papers to study and constrain transport theories such as diffusion phenomena across or along resonances \cite{cFr00,mGu13,eLe16}.  The problem is trivially integrable when $\epsilon=0$. For $\epsilon \neq 0$, 
the unperturbed resonances associated to the Hamiltonian of Eq.\,(\ref{eq:FGL}) read 
\begin{align}\label{eq:ResFGL}
k \cdot \partial_{I} \mathcal{H}_{0}(I,\phi) =
k_{1}I_{1} + k_{2}I_{2} + k_{3}, \, k=(k_{1},k_{2},k_{3}) \in \mathbb{Z}_{\star},
\end{align}
and translate as straight lines in the $(I_{1},I_{2})$ action space. The set of resonances is dense within this plane but their effects decrease with the order of the resonance 
$\vert k \vert$. To reveal the resonant template, the interactions among the resonances and its parametric evolution according to $\epsilon$,  we perform a $\norm{\Delta \LD}$ stability 
analysis  following strictly the seminal work of \cite{cFr00}. In the latter,  the FLI is used to portray the evolution of the resonant  web. 
Fig.\,\ref{fig:fig3} shows the evolution of the resonant web for increasing values of the perturbing parameter ($\epsilon=0.001$, $\epsilon=0.01$ and $\epsilon=0.04$) at different scales of the action space. 
The macroscopic domain is defined by the section
\begin{align}
\Sigma_{\textrm{M}}=\big\{ (I_{1},I_{2}) 
\in [-0.5, 1.5]^{2},\,
I_{3}=\phi_{1}=\phi_{2}=\phi_{3}=0
\big\},
\end{align}
and corresponds to the left panel of Fig.\,\ref{fig:fig4}. Enlargements of this section define the microscopic section defined as 
\begin{align}
\Sigma_{\textrm{m}}=
\big\{ (I_{1},I_{2}) 
\in [0.3, 0.4] \times [0.1, 0.2],\,
I_{3}=\phi_{1}=\phi_{2}=\phi_{3}=0
\big\}
.
\end{align}
The dynamical portrays of this scale correspond to the right part of Fig.\,\ref{fig:fig3}. 
The parametric $\norm{\Delta \LD}$ analysis reproduces in detail the result of \cite{cFr00}. For the small value  $\epsilon=0.001$ (top row of Fig.\,\ref{fig:fig3}), the phase space is predominantly filled by regular motions. One detects the presence of many resonances corresponding to the lines of Eq.\,(\ref{eq:ResFGL}). The volume of regular orbit decreases for larger $\epsilon$, as made evident for $\epsilon=0.01$ (middle row of Fig.\,\ref{fig:fig3}). Chaotic motions appear sharply at low-order resonant crossings. This is especially visible at the microscopic scale, where also many thin secondary substructures are detected. The volume of regular orbits shrinks further for larger value of $\epsilon$. 
At $\epsilon=0.04$ (bottom row of Fig.\,\ref{fig:fig3}), chaotic motions are also found at higher-order resonant crossings. 
At the lowest scale, one notices the presence of larger chaotic seas allowing faster transport routes in the phase space as the result of resonances overlap. 
This dynamical regime is substantially different from the regime where resonances are well separated. 
The detection of chaotic and regular orbits, and its spatial arrangement, is trustingly recovered by $\norm{\Delta \LD}$.    

\section{Applications to mappings}\label{sec:Map}
This section demonstrates the ability of the $\norm{\Delta \LD}$ indicator to reveal accurately phase space structures of discrete systems. We utilise symplectic and volume-preserving mappings as test beds.   

\subsection{The standard map}\label{subsec:SM}
\begin{figure}
	\centering
	\includegraphics[width=1\linewidth]{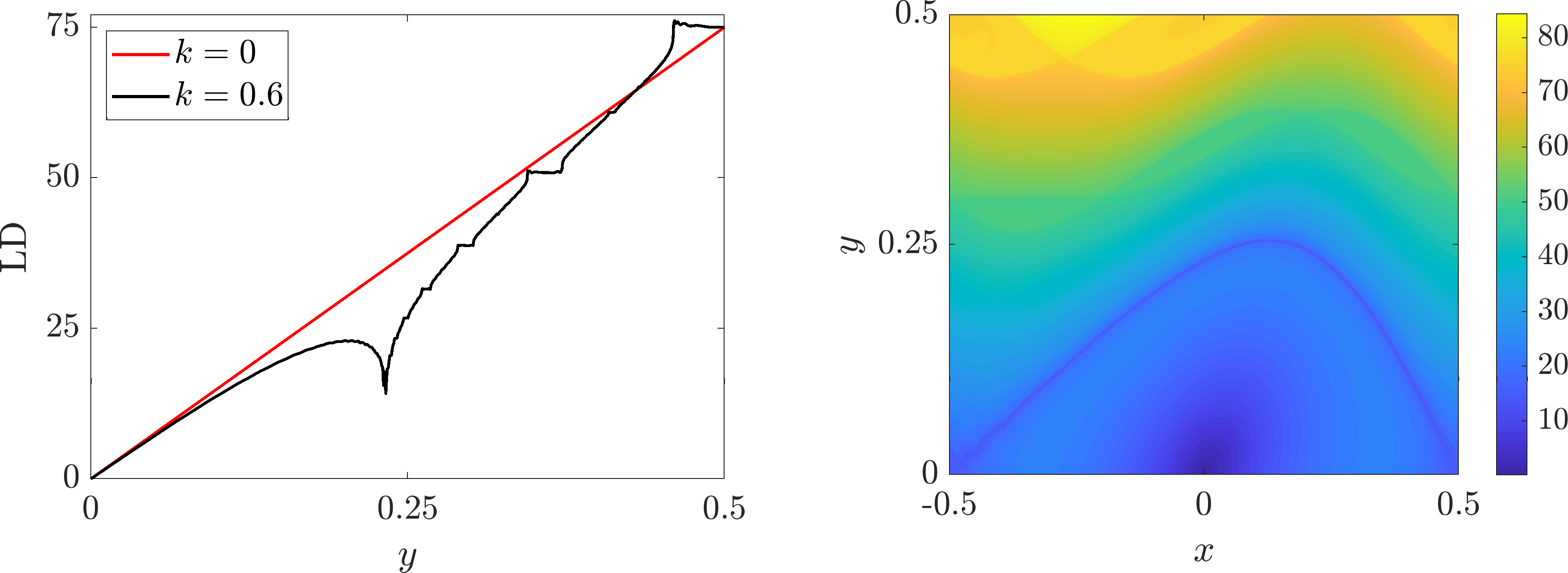}
	\caption{
		(Left)  LD landscape of Eq.\,(\ref{eq:SM}) for $k=0$ (integrable case) and $k=0.6$ computed at $n=150$ for $x=0$.
		(Right) LD field of Eq.\,(\ref{eq:SM}) for $k=0.6$ computed at $n=150$. 
		Albeit the information about the locations of the hyperbolic structures is contained within the LD metric (loss of regularity), the heatmap of the LD field itself is not able to restore them sharply.
		The various separatrices are drowning in the linear trend of the LD metric as a function of $y_{0}$ (confer Eq.\,(\ref{eq:LDint})), thus contributing to a ``flat'' map effect.}
	\label{fig:fig4}
\end{figure}

\begin{figure}
	\centering
	\includegraphics[width=1\linewidth]{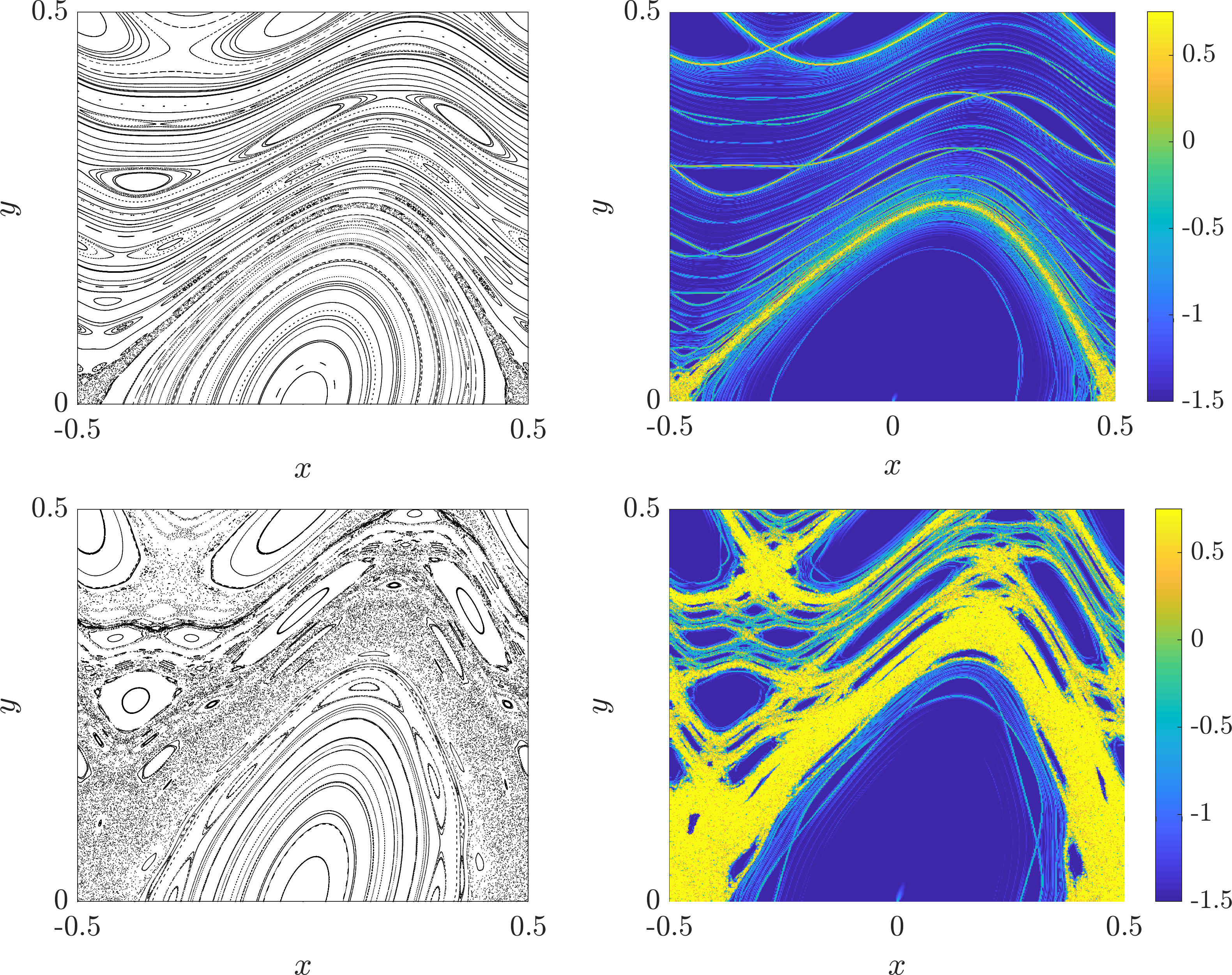}
	\caption{
		(Top line) Phase space of the standard map for $k=0.6$ obtained by 
		iterating trajectories up to the final time $n=750$ (left panel) 
		and using the $\norm{\Delta \LD}$ indicator computed at $n=150$ (right panel).
		(Bottom line) Same computations performed with $k=1$.  	
	}\label{fig:fig5}
\end{figure}

The paradigmatic standard map is defined on $\mathbb{T}\times\mathbb{R}$ as
\begin{align}
(x,y)\mapsto (x',y')=f_{k}(x,y), 
\end{align}
with
\begin{align}\label{eq:SM}
f_{k}:
\left\{
\begin{aligned}
& x'= x+ y + F_{k}(x) \mod 1,\\
& y'= y + F_{k}(x),
\end{aligned}
\right.
\end{align}
where $F_{k}(x)=-k\sin(2\pi x)/(2\pi)$, 
$k \in \mathbb{R}^{+}$ is the nonlinearity parameter. 
For $k=0$, the map $f_{0}$ becomes  
\begin{align}\label{eq:SMIntegrable}
f_{k}:
\left\{
\begin{aligned}
& x'= x+ y  \mod 1,\\
& y'= y,
\end{aligned}
\right.
\end{align}
and is integrable. 
The solutions at time $n$, starting from the initial condition $(x_{0},y_{0})$, read
\begin{align}\label{eq:SMint}
x_{n}=x_{0} + n y_{0}, \, \, y_{n}=y_{0}, 
\end{align} 
leading to 
\begin{align}\label{eq:LDint}
\LD\big((x_{0},y_{0}),n\big) =
\sum_{i=0}^{n-1}
\sqrt{
	(x_{i+1}-x_{i})^{2} + (y_{i+1}-y_{i})^{2}
}
=
n \sqrt{y_{0}^{2}}
=
n \vert y_{0} \vert.
\end{align}  
When $k \neq 0$, the phase space contains a mixture of invariant curves and chaotic motions, densely filling the phase space as $k$ increases \cite{jMe92,jMe08}. \nico{This map using a different version of LD has been also discussed by \cite{Bib3}.}
Eq.\,(\ref{eq:LDint}) represents the discrete analogue of Eq.\,(\ref{eq:LDlinearInte}) and  encapsulates the inherent limitations of portraying the LD field through a heatmap to visualise the dynamical structures. The left panel of Fig.\,\ref{fig:fig4} shows LD landscapes computed for $x=0$, up to the time $n=150$, for the integrable $k=0$ case  (red line) and $k=0.6$ (black line) for $y \in \mathcal{D}=[0,0.5]$. As predicted by Eq.\,(\ref{eq:LDint}), the LD of the integrable case grows linearly as a function of the initial action $y$.  
Increasing $k=0$ to $k=0.6$, we observe that the LD landscape  is predominantly guided by the integrable approximation.
The most noticeable difference to this trend occurs in the vicinity of the main resonant island and its ``separatrix.'' 
Although the location of the manifolds is guessable by visual inspection of the LD landscape, the small amplitudes variations of the LDs when crossing the hyperbolic layers combined with the overall sharp linear trend of the LDs  tend to erase the information (presence of two distinct scales). Consequently, the heatmap of the LD field is ``flat'' and unstructured, in the sense that it does not contain clear ridges
associated to the dynamical structures, as shown in the right panel of Fig.\,\ref{fig:fig5}. \\

Fig.\,\ref{fig:fig5} compares the phase space analysis using iterations of orbits and demonstrates that the $\norm{\Delta \LD}$ indicator alleviates this issue (compare the top right plot of Fig.\,\ref{fig:fig5} with the right panel of Fig.\,(\ref{fig:fig4})). 
The phase spaces have been obtained by iterating $f_{k}$  for $n=750$ times. 
For $k=0.6$, the phase space contains primarily invariant librational and circulational curves. 
For $k=1$, resonant domains have expanded and overlapped, and are no longer separated by invariant curves. 
The phase space contains a larger volume of chaotic orbits. 
The corresponding results of the  $\norm{\Delta \LD}$ analysis  is shown  in the right part of Fig.\,\ref{fig:fig5}. 
The  $\norm{\Delta \LD}$ are computed for the final time $n=150$ on a regular mesh of $500 \times 500$ initial conditions.  
Indisputably, the heatmap of the $\norm{\Delta \LD}$ indicator succeeds in recovering the global template of the system. The topology of the lobes of the various resonances, and the distribution of chaos around them, are clearly revealed.   

\subsection{A $4$-dimensional nearly-integrable mapping}\label{subsec:4d}
We now turn our attention to a higher dimensional discrete system, namely the  symplectic nearly-integrable mapping defined by:
\begin{align}\label{eq:MFroeschle}
\left\{
\begin{aligned}
&x_{j+1} = x_{j} -\epsilon \sin(x_{j}+y_{j})/\mu(x_{j},y_{j},z_{j},t_{j}), \\
&y_{j+1} = y_{j} + x_{j}, \\
&z_{j+1} = z_{j} - \epsilon \sin(z_{t}+t_{j})/\mu(x_{j},y_{j},z_{j},t_{j}), \\
&t_{j+1} = z_{j} + t_{j},
\end{aligned}
\right.
\end{align}
with
\begin{align}
\mu(x_{j},y_{j},z_{j},t_{j})=\big(\cos(x_{j}+y_{j}) + \cos(z_{j}+t_{j}) + 4\big)^{2}. 
\end{align}
When $\epsilon=0$, the mapping is integrable. Similarly to the standard map, $x$ and $z$ are both constant, whilst $y$ and $t$ evolve linearly with time. 
This mapping, and variations of it, have been classical molds to study transport and diffusion phenomena along resonances in nearly integrable settings \cite{cFr05,mGu06}.  
The resonances associated to the system of Eq.\,(\ref{eq:MFroeschle}) read \cite{mGu04}
\begin{align}
k_{1}x + k_{2}z + 2 k_{0} \pi = 0, \, (k_{1},k_{2},k_{0}) \in \mathbb{Z}_{\star}^{3},
\end{align}
and translate as straight lines into the $x-z$ plane. 
Albeit the set of resonances is dense into this plane, resonant orbits surround the resonant locations with a distance that decreases with the order of the resonance $\vert k \vert$. Understanding analytically the multi-resonant dynamics, \ie the locations of the resonances, the strength of each of them and their mutual interactions in a hierarchical way, is a difficult task \cite{aMo95}. Instead, to reveal the hyperbolic structures  one follows here a purely numerical procedure by exploring the \textit{geography of resonances} using the $\norm{\Delta \LD}$ indicator.  \\

\begin{figure}
	\centering
	\includegraphics[width=1\linewidth]{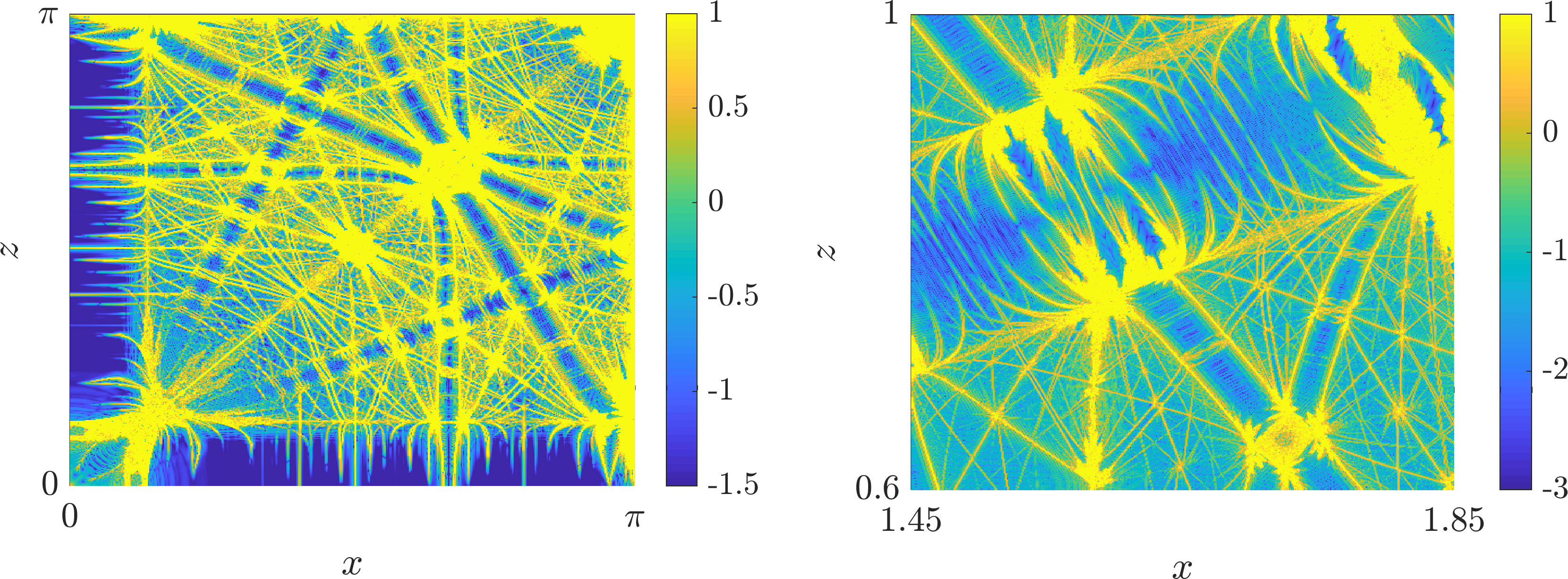}
	\caption{
		$\norm{\Delta \LD}$ maps associated to the mapping of Eq.\,(\ref{eq:MFroeschle}) for $\epsilon=0.6$ computed 
		on $\Sigma_{1}$ and $\Sigma_{2}$ at $t=1,000$.
		The $\norm{\Delta \LD}$ indicator recovers the resonant templates and the presence of hyperbolic orbits.	
	}\label{fig:fig6}
\end{figure}

Fig.\,\ref{fig:fig6} presents the details of the geography of the resonances associated to Eq.\,(\ref{eq:MFroeschle}) with $\epsilon=0.6$ computed at $n=1,000$ at two different scales of the phase space. The numerical settings follow \cite{cFr05}. The sections are respectively defined by
\begin{align}
\Sigma_{1} = \big\{(x,y,z,t) \, \vert \, (x,z) \in [0,\pi]^{2}, \, y=t=0\big\}, 
\end{align}
and 
\begin{align}
\Sigma_{2} = \big\{(x,y,z,t) \, \vert \, (x,z) \in [1.45,1.85]\times [0.6,1], \, y=t=0\big\}. 
\end{align}
The reduced scale $\Sigma_{2}$ focuses on the resonant structure along the $x=2z$ resonance. 
The results of this $\norm{\Delta \LD}$ analysis are in excellent agreement with the FLI maps produced in the Fig.\,2 and Fig.\,3 of \cite{cFr05}. 
The $\norm{\Delta \LD}$ indicator reflects at the macroscopic scale the predominance  of  chaos near low-order resonance crossings.  
At the lower scale, $\norm{\Delta \LD}$ is able to detect sharply details of the rich dynamical structure, especially the ``background'' of high-order resonances. 

\subsection{Froeschl\'e's generalised map}\label{sub:GFM}
\begin{figure}
	\centering
	\includegraphics[width=0.8\linewidth]{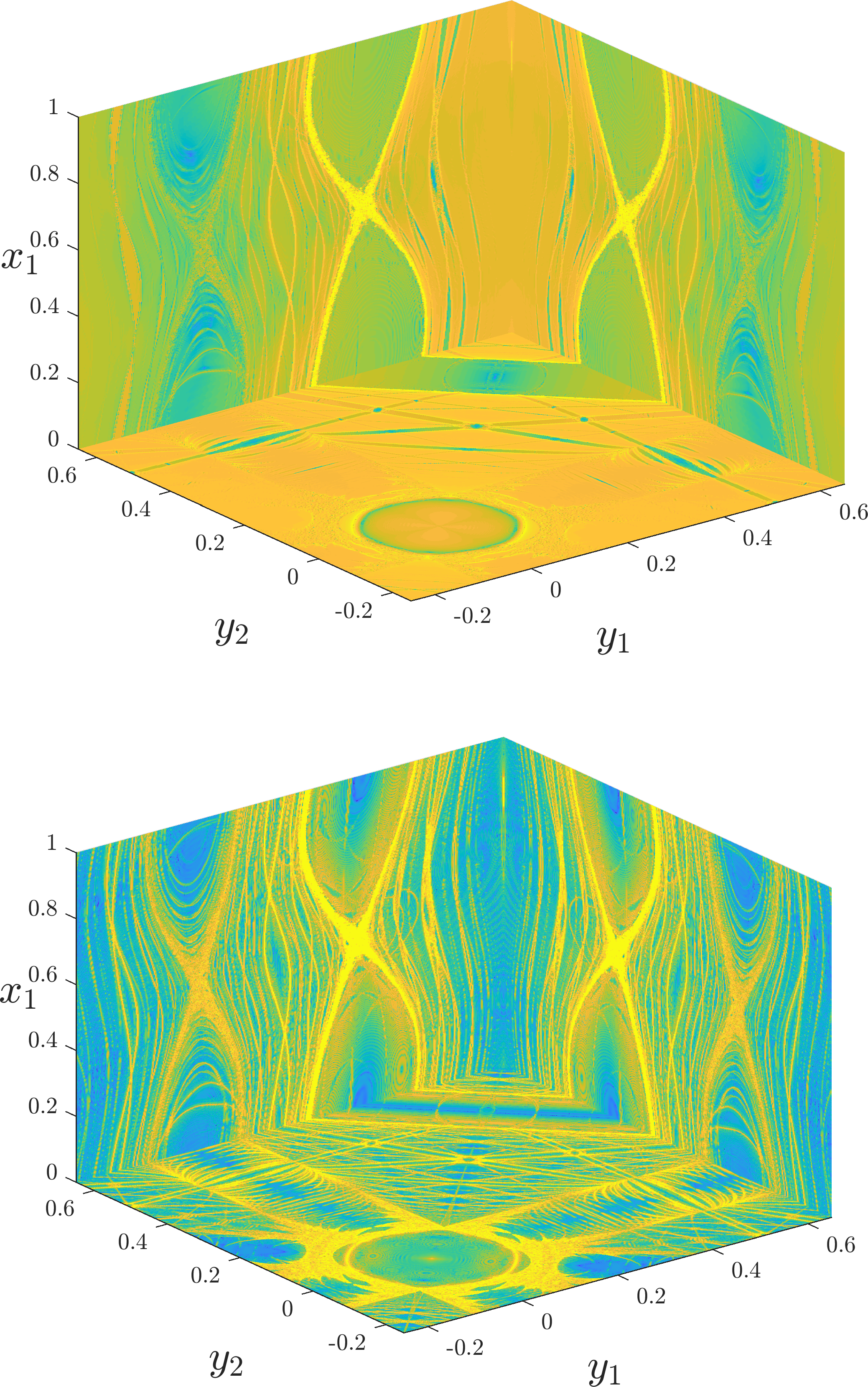}
	\caption{
		Dynamical maps associated to Eq.\,(\ref{eq:FM}) in the symplectic case ($\varphi=0$) using
		(top) $\norm{\nabla \LD}$ and (bottom) $\norm{\Delta \LD}$. The second-order derivatives based indicator is better suited than the gradient approach to portray the geography of resonances. 
	}
	\label{fig:fig7}
\end{figure}

\begin{figure}
	\centering
	\includegraphics[width=0.8\linewidth]{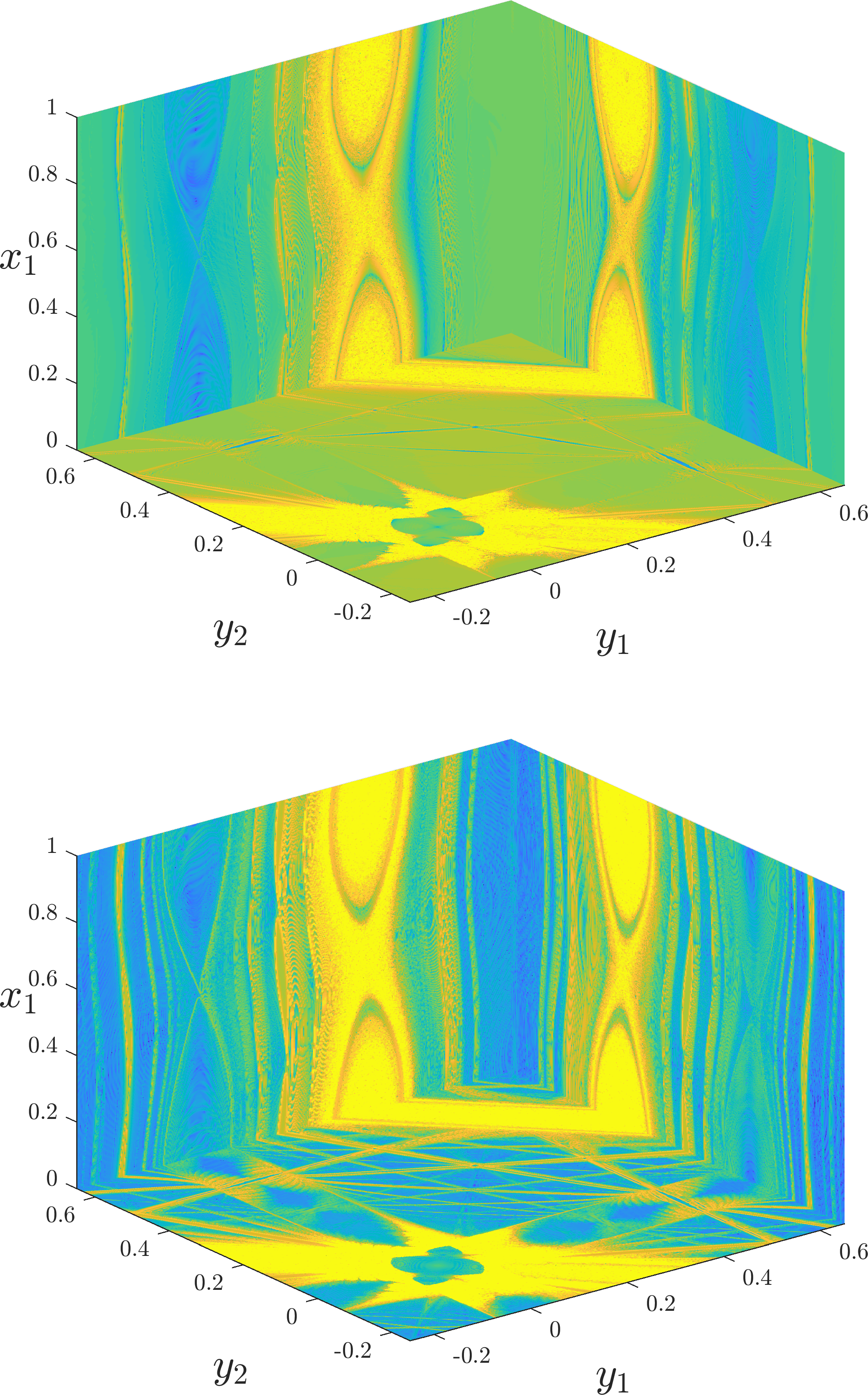}
	\caption{
		Dynamical maps associated to Eq.\,(\ref{eq:FM}) in the volume-preserving case ($\varphi=0.1$) using
		(top) $\norm{\nabla \LD}$ and (bottom) $\norm{\Delta \LD}$. The second-order derivatives based indicator is better suited than the gradient approach to portray the geography of resonances.
	}
	\label{fig:fig8}
\end{figure}

The $4$-dimensional generalised Froeschl\'e's map defined on 
$\mathbb{T}^{2}\times\mathbb{R}^{2}$ by 
\begin{align}\label{eq:FM}
f_{(a,b,c)}:
\left\{
\begin{aligned}
& x'_{1}= x_{1}+ y_{1} - \frac{1}{2\pi} 
\Big(
a \sin(2\pi x_{1}) +
c \sin(2\pi(x_{1}+x_{2})) \Big)
\mod 1,\\
& x'_{2}= x_{2}+ y_{2} - \frac{1}{2\pi} 
\Big(b \sin(2\pi x_{2}) +
c \sin(2\pi(x_{1}+x_{2})) \Big)\mod 1,\\
& y'_{1}= y_{1} - \frac{1}{2\pi} 
\Big(a \sin(2\pi x_{1}) +
c \sin(2\pi(x_{1}+x_{2}))\Big),\\
& y'_{2}= y_{2} - \frac{1}{2\pi}
\Big(
b \sin(2\pi x_{2}) +
c \sin(2\pi(x_{1}+x_{2}+\varphi))\Big),
\end{aligned}
\right.
\end{align}
has been proposed by \cite{nGu17} to compare the dynamics and transport properties of symplectic and volume preserving maps. 
The parameters $(a,b,c)$ represent forcing terms of the $(1,0,n)$, $(0,1,n)$  and $(1,1,n)$ resonances respectively, $n \in \mathbb{Z}$. 
The map is symplectic for $\varphi=0 \mod 1$, and volume-preserving for nonzero $\varphi$. 
The parameter $\varphi$ is thus a convenient measure of deviation from symplecticity.  
When $c=0$, Eq.\,(\ref{eq:FM}) becomes decoupled standard maps, as discussed in Sect.\,\ref{subsec:SM}. 
Fig.\,\ref{fig:fig7} and \ref{fig:fig8} are reproduction of the visualistion of the high-dimensional dynamics, following the steps of 
\cite{nGu17}, comparing the first-order derivatives based $\norm{\nabla \LD}$ and $\norm{\Delta \LD}$. 
The presentation of the maps into a cube is a convenient condensation to appreciate the need of the second-order based diagnostic to delineate the geography and interactions among resonances. 
Whilst $\norm{\Delta \LD}$ succeeds in recovering the dynamical template, $\norm{\nabla \LD}$ misses the main and secondary resonant strips in the action plane. The situation improves in the angle-action space, yet small resonant islands are undetected.  
The parameters of the simulations read $(a,b,c)=(0.1,0.1,0.07)$, $\varphi=0$ and $(a,b,c)=(0.05,0.05,0.035)$, $\varphi=0.1$ respectively.   
Each map  is the result of $500 \times 500$ initial conditions, propagated up to the final time $n=1,000$.

\section{Summary and conclusive remarks}
This paper has introduced a new global dynamics  and chaos indicator based on the theory of Lagrangian Descriptors. 
The $\norm{\Delta \LD}$ indicator, based on the second-derivatives of the LDs, \nico{decapsulates the lack of regularity of the LD metric allowing a visualisation of global dynamical features}. The main contributions and conclusions  are the following:
\begin{enumerate}
	\item Producing a heatmap based on the LD field itself might not reveal all the information about the locations and the precise geometry of hyperbolic structures (a phenomena described as ``flat maps.'') 
	This is even more true when several nearby hyperbolic structures cohabit in the phase space. In the nearly-integrable setting we investigated here, this fact is explained by the linear growth of the LD with respect to the initial ``action,'' which tends to overwhelm relevant dynamical information. 
	The proposed $\norm{\Delta \LD}$ indicator  alleviates this limitation.
	\item The $\norm{\Delta \LD}$ indicator is apt to discriminate between ordered and chaotic motions. The indicator has been benchmarked on several multidimensional continuous and discrete  models against several phase space methods. 
	These included the computation of dynamical maps of various perturbed pendulums, 
	the $2$-\DOF \,H\'enon-Heiles system and a $3$-\DOF \, Hamiltonian system supporting a dense web of resonances and diffusive phenomena. The diagnostic has been validated for the discrete realm, by computing stability maps of the 
	standard map and a $4$ dimensional nearly-integrable mapping.
	Our set of $\norm{\Delta \LD}$ maps can be confronted to existing maps produced with various variational methods, such as the FLI, the MEGNO, the orthogonal FLI, or traditional phase space techniques such as the
	phase space method or iterations of the period-map. 
	Our simulations demonstrate in particular that resonant and chaotic templates can be recovered through $\norm{\Delta \LD}$ cartography. The $\norm{\Delta \LD}$ method recovers minutes details of the dynamics across several time and space scales, and is successful in delineating the geography of resonances.  
\end{enumerate}

The $\norm{\Delta \LD}$ indicator, derived from the LD metric, does not rely on the concept of separation of nearby orbits and its quantification through the growth of the norm of the tangent vector. Its implementation is thus free of the variational equations, and requires only to compute arc-lengths of trajectories  on calibrated finite size time-windows. This property, convenient by itself (only the level of the equation of motions is needed), also implies a computational advantage over variational methods \nico{as already reported in \cite{aMa13,aDa21}. A more precise quantification will be provided in a forthcoming paper.} 
The newly introduced $\norm{\Delta \LD}$ indicator is able to unveil resonant and chaotic templates, and thus 
appears to be  relevant for the fields of celestial mechanics, dynamical astronomy and astrophysics for studying problems related to transport in the phase space shaped by resonant interactions.     

\section*{Acknowledgments}
\nico{The authors are very grateful to  V\'ictor J. Garc\'ia-Garrido and Stephen Wiggins for bringing to their knowledge their recent references \cite{vGa22,mKa22} and useful comments and suggestions}.  
J.\,D. acknowledges warmly discussions and feedback from Carolina Charalambous, Anne Lemaitre and Timoteo Carletti. 
J.\,D. is a postdoctoral researcher of the ``Fonds de la Recherche Scientifique'' - FNRS.
M.\,A. acknowledges support from the grant CEX2019-000904-S and IJC2019-040168-I funded by: MCIN/AEI/ 10.13039/501100011033 \nico{by “European Union NextGenerationEU/PRTR”}.
\nico{A.M.\,M acknowledges  support from grant PID2021-123348OB-I00 funded by  
	MCIN/ AEI /10.13039/501100011033/ and by
	FEDER A way to make Europe.}
\appendix

\section{Calibration of the time window \& geometrical LDs}\label{app:LDgeo}
In the first line of Fig.\,\ref{fig:fig1} is computed arc-lengths of orbits for the integrable pendulum model $\mathcal{H}$ over the time window $[0,t]$,  $t=100$.  
Alike many chaos detection methods, there is no strict theoretical guidance for the choice of the final time $t$. 
The practitioner might take advantage of the knowledge of 
some specific timescale,  or perform others simulations and saturation checks in order to calibrate this time window. 
A too short time misses the detection of the structures, whilst a prohibitive large time increases the computational burden.   
To get rid of this time dependence, we have developed a geometrical framework of the Lagrangian Descriptor, called \textit{geometrical Lagrangian Descriptor}, for the class of $1$-\DOF \, Hamiltonian system \cite{rPO21}. 
This framework and point of view has several benefits. Firstly, the  lengths are no longer parameterized by the time but only by the energy of the orbit of the system (the geometrical LDs are thus completely free of the time variable). Given an energy level $E$, the geometrical LD associated to $E$, denoted $\ell(E)$, corresponds to the length of the level curve $\mathcal{H}(I,\phi)=E$. 
In the case of the pendulum model,  the level curves  $\mathcal{H}(I,\phi)=E$ on $\mathbb{R} \times [-\pi,\pi]$ are interpreted as the planar curves parametrised by $\phi$ given by 
\begin{align}
I(\phi;E)=\pm \sqrt{2(E+\cos \phi)}. 
\end{align} 
Exploiting the formula for the length of a curve and symmetries in the phase space, we might rewrite $\ell(E)=2\tilde{\ell}(E)$
with
\begin{align}\label{eq:lE}
\tilde{\ell}(E)= 
\int 
\sqrt{1 + \Big( 
	\frac{\rm{d} \textit{I}}{\rm{d} \phi}\Big)^{2}} \, \rm{d} \phi,
\end{align} 
where the integral is computed over a suitable range for $\phi$. 
Formula in Eq.\,(\ref{eq:lE}) has a deeper geometrical content  than Eq.\,(\ref{eq:LD}). In particular, it depends solely on the energy of the system. The time-free analogue of the landscape presented in the top right panel of Fig.\,\ref{fig:fig1}, 
$\ell(E)$ as a function of $E$, has also a cusp point at the energy $E$ labeling the separatrix (note that in the example we investigated, we found that the separatrix curve always realises a local maximal of $\ell(E)$). 
We refer to  \cite{rPO21} for omitted details or complementary figures related to the pendulum model and other integrable problems on the cylinder and plane (including the $8$ and fish-shaped separatrices). Secondly, the geometrical apparatus has allowed to characterise, analytically or semi-analytically, the speed at which $\vert \dd \ell(E)/\dd E \vert$ becomes singular when
$E \to E_{\textrm{sx}}$, where  $E_{\textrm{sx}}$ denotes the energy level labeling the separatrix. For the examples considered and independently of the separatrices topology,  we always found power-laws scaling as $\mathcal{O}(1/\sqrt{\vert E-E_{\textrm{sx}}\vert})$.

\begin{remark}[LD based on the actions] 
	In the setting of a $n$-\DOF Hamiltonian having the form ``kinetic + potential'' energy,
	\begin{align}
	\mathcal{H}(p,q)=T(p)+V(q),
	\end{align}  
	where $T(p)=\sum_{i=1}^{n}p_{i}^{2}/2$,
	the recent work of \cite{gGa22} proposed the reduced Maupertuis action
	\begin{align}\label{eq:LDaction}
	\LD(x_{0},t)=\int_{q_{0}}^{q_{t}} p \,\dd q,
	\end{align}
	as Lagrangian Descriptor. Here $x_{0}=(p_{0},q_{0})$ is the initial condition at time $s=0$, $q_{0}$ and $q_{t}$ refer to the configuration 
	coordinates at the initial time $s=0$ and final time $s=t$ of the time window $[0,t]$.
	Whilst Eq.\,(\ref{eq:LDaction}) resembles the geometrical LD  introduced by \cite{rPO21}, the authors then exploit Hamilton's canonical equations to rewrite  Eq.\,(\ref{eq:LDaction}) as
	\begin{align}
	\LD(x_{0},t)=\int_{q_{0}}^{q_{t}} p \,\dd q
	=
	\int_{0}^{t} p\,\frac{\dd q}{\dd s}\,\dd s
	=
	\int_{0}^{t} 2\,T(s)\,\dd s,
	\end{align}
	which ends up being a LD computed ``classically'' in the time domain. 
\end{remark}

\section{Application to the H\'enon-Heiles system}\label{app:HH}
The $2$-\DOF \,H\'enon-Heiles system is a paradigmatic example of Hamiltonian chaos that has received substantial theoretical and computational attention. The H\'enon-Heiles Hamiltonian, rooted in galactic dynamics, reads \cite{mHe64}
\begin{align}\label{eq:HHH}
\mathcal{H}(x,y,p_{x},p_{y})
=
\frac{1}{2}(p_{x}^{2}+p_{y}^{2}+x^{2}+y^{2}) + x^{2}y - \frac{1}{3}y^{3}. 
\end{align}

\begin{figure}
	\centering
	\includegraphics[width=1\linewidth]{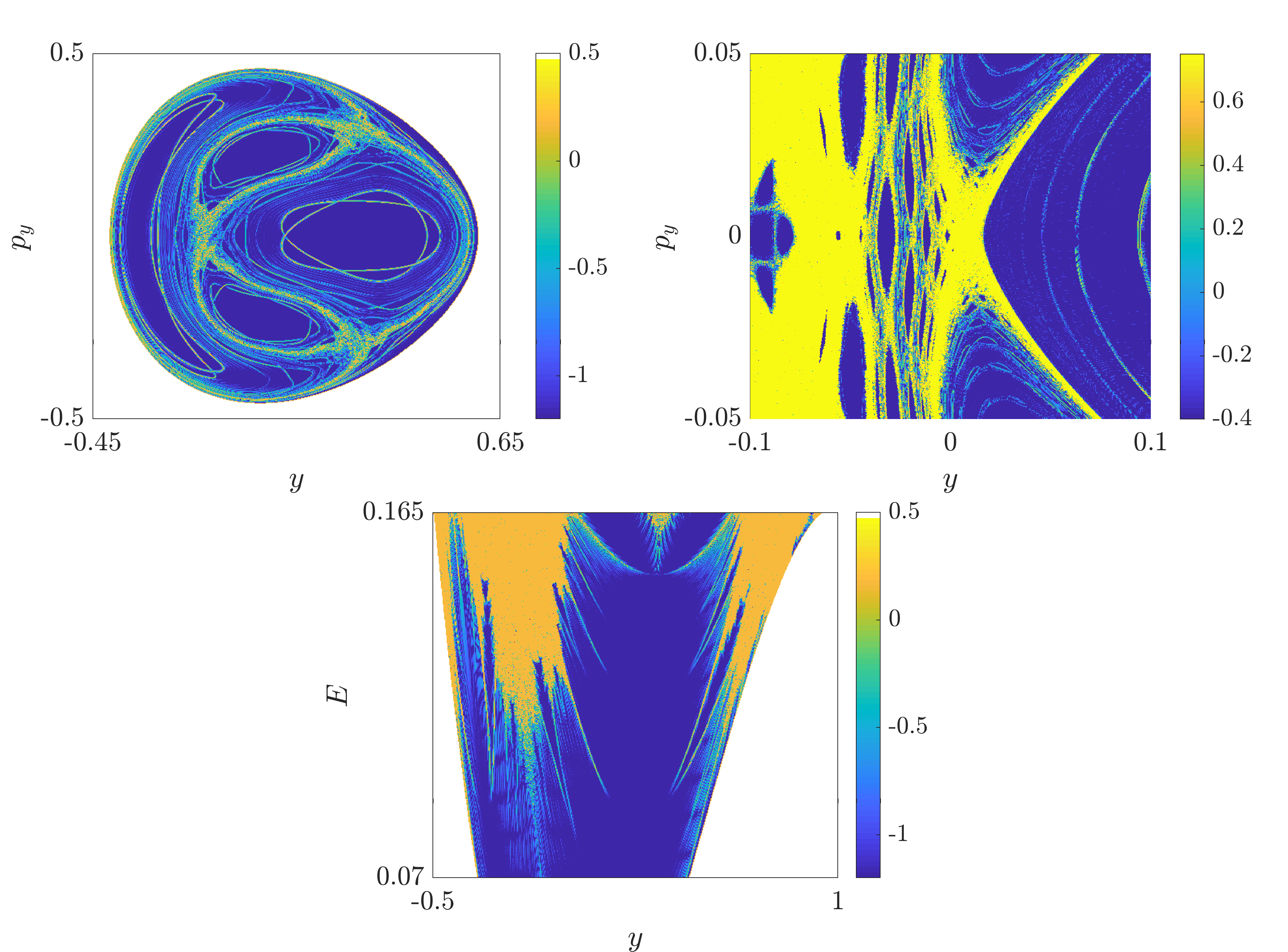}
	\caption{
		Dynamical maps associated to the H\'enon-Heiles system given in Eq.\,(\ref{eq:HHH}) computed in various planes with the $\norm{\Delta \LD}$ indicator. 
		The $\norm{\Delta \LD}$ indicator succeeds in revealing the dynamical structures at several spatio and temporal scales.  
	}
	\label{fig:fig9}
\end{figure}

Fig.\,\ref{fig:fig9} reproduces $3$ stability maps associated to the system (\ref{eq:HHH}) at different scales following strictly the numerical setting of  \cite{rBa05,pCi16}. In the latter, the orthogonal FLI and the MEGNO indicator are the variational tools respectively used to portray the phase spaces. Before commenting further the results, let us emphasise that in all the heatmaps of panel \nico{\ref{fig:fig9},} a white color is assigned to non-admissible couple $(y,p_{y})$ (this will be clearer in the subsequent). \nico{This system has been also discussed with different implementations of LDs in \cite{Bib4, Bib1,Bib3}}. 
We summarise the results of the numerical procedure as follow:
\begin{enumerate}
	\item The top-left map of Fig.\,\ref{fig:fig9} presents the global phase space at the energy value $E=0.105$ following \cite{rBa05} where the orthogonal FLI is used to portray the phase space. The final time $t$ is set to $t=300$ and we use a $500 \times 500$ grid of initial conditions. 
	The section on which the $\norm{\Delta \LD}$ is computed 
	let free the variables $(y,p_{y})$ with  $x=0$. The last variable $p_{x}$ to initialise the differential system is determined by solving the iso-energetic equation $\mathcal{H}(x,y,p_{x},p_{y})=E$,   $p_{x} > 0$. If the latter equation does not admit a solution, the uplet $(y,p_{y})$ is said to be non-admissible. The  
	$\norm{\Delta \LD}$ map reproduces faithfully the result of \cite{rBa05}. 
	\item The top-right map of Fig.\,\ref{fig:fig9} focuses on a smaller portion of the previous section (\ie $(y,p_{y})$, $x=0$ and $p_{x} > 0$ is determined by the iso-energetic condition) for a slightly higher value of energy, $E=0.118$. The numerical setting follows \cite{pCi16} for which the MEGNO indicator was computed.   
	The final time of the computation of the $\norm{\Delta \LD}$  map is $t=10^{4}$, \ie about two orders of magnitude larger than  the former chosen time. 
	The map thus offers a resolved view of the long-term dynamics near the origin of the system.    
	The result of the $\norm{\Delta \LD}$ map demonstrates the ability of the LD metric to recover tiny structures at a very fine scale, as observed from the minute lobes that are distinguishable in the vicinity of the origin.
	The $\norm{\Delta \LD}$ indicator offers a clear picture of the dynamics, and is able on this example to deliver more details than the MEGNO analysis (both are computed for the same final time), in particular, in detecting the thin unstable domain within the stable island.    
	\item The last map  of Fig.\,\ref{fig:fig9} presents stability results in the complimentary section $(y,E)$. As before, we set $x=0$, $p_{y}=0$, and $p_{x} > 0$ is determined by the energy condition. The final time is set to $t=300$. 
	The result of the $\norm{\Delta \LD}$ computation, that can be compared to  \cite{rBa05}, is in excellent agreement all along the energy range probed. In particular, the analysis is able to recover the tiny fringes of instability (tongues) foliating the stable domain.   
\end{enumerate}

\bibliographystyle{plain}   
\bibliography{biblio}

\end{document}